\documentclass[amssymb,amstex,amsmath,10pt,a4paper]{amsart}
\usepackage{amsmath,amssymb,latexsym,amsfonts,url}

\newtheorem{thm}{Theorem}[section]
\newtheorem*{theorem}{Theorem}

\newtheorem{cor}[thm]{Corollary}

\newcommand{\area}{{\rm Area}}
\newcommand{\length}{{\rm Length}}

\newcommand{\Ric}{{\rm Ric}}
\newcommand{\Div}{{\rm div}}
\newcommand{\ind}{{\rm Index}}

\newenvironment{pf}{\noindent {\it Proof.}}{\\}
\newenvironment{Rmk}{\noindent {\it Remark.}}{\\}
\numberwithin{equation}{section}

\begin{document}
\title[{Stable minimal hypersurfaces in the hyperbolic space}]{Stable minimal hypersurfaces \\in the hyperbolic space}
\author{Keomkyo Seo} \address{Department of Mathematics\\
Sookmyung Women's University\\
Hyochangwongil 52, Yongsan-ku \\
Seoul, 140-742, Korea}
\email{kseo@sookmyung.ac.kr}


\begin{abstract}
In this paper we give an upper bound of the first eigenvalue of
the Laplace operator on a complete stable minimal hypersurface $M$
in the hyperbolic space which has finite $L^2$-norm of the second
fundamental form on $M$. We provide some sufficient conditions for
minimal hypersurface of the hyperbolic space to be stable. We also
describe stability of catenoids and helicoids in the hyperbolic
space. In particular, it is shown that there exists a family
of stable higher-dimensional catenoids in the hyperbolic space. \\

\end{abstract}

\maketitle

\section{Introduction}
In \cite{CLY}, Cheng, Li and Yau derived comparison theorems for
the first eigenvalue of Dirichlet boundary problem on any compact
domain in minimal submanifolds of the hyperbolic space by
estimating the heat kernel of the compact domain. Recall that the
first eigenvalue $\lambda_1$ of a complete non-compact Riemannian
manifold $M$ is defined by $\lambda_1 = \inf_{\Omega}
\lambda_1(\Omega)$, where the infimum is taken over all compact domains in $M$. Throughout this paper, we shall denote by
$\mathbb{H}^n$ the $n$-dimensional hyperbolic space of constant
sectional curvature $-1$. Recently Candel \cite{Candel} gave un
upper bound for the first eigenvalue of the universal cover of a
complete stable minimal surface in $\mathbb{H}^3$. Indeed, he
proved
\begin{theorem}[\cite{Candel}]
Let $\Sigma$ be a complete simply connected stable minimal surface
in the $3$-dimensional hyperbolic space. Then the first eigenvalue
of $\Sigma$ satisfies
\begin{eqnarray*}
\frac{1}{4} \leq \lambda_1 (\Sigma) \leq \frac{4}{3} .
\end{eqnarray*}
\end{theorem}

In Section 2, we extend this theorem to simply connected stable
minimal surfaces in a Riemannian manifold whose sectional
curvature is bounded below and above by negative constants
(Theorem \ref{thm:manifold}). For a complete stable minimal
hypersurface $M$ in $\mathbb{H}^{n+1}$, Cheung and Leung
\cite{Cheung and Leung} proved that
\begin{eqnarray} \label{thm : CL}
\frac{1}{4}(n-1)^2 \leq \lambda_1 (M) .
\end{eqnarray}
Here this inequality is sharp because equality holds when $M$ is
totally geodesic (\cite{McKean}). In this paper, it is proved that if $M$ is a
complete stable minimal hypersurface in $\mathbb{H}^{n+1}$ with
finite $L^2$-norm of the second fundamental form $A$, then we have
(Theorem \ref{thm:first eigenvalue in H})
\begin{eqnarray*}
\lambda_1(M) \leq n^2 .
\end{eqnarray*}
Recall that a minimal hypersurface is called {\it stable} if the second
variation of its volume is always nonnegative for any normal
variation with compact support. More precisely, an $n$-dimensional
minimal hypersurface $M$ in a Riemannian manifold $N$ is called
{\it stable} if it holds that for any compactly supported
Lipschitz function $f$ on $M$
\begin{eqnarray} \label{def : stability}
\int_M |\nabla f|^2 - \Big(|A|^2 + \overline{\Ric}(\nu,\nu)\Big) f^2
dv \geq 0 ,
\end{eqnarray}
where $\nu$ is the unit normal vector of $M$,
$\overline{\Ric}(\nu,\nu)$ denotes the Ricci curvature of $N$ in
the $\nu$ direction, $|A|^2$ is the square length of the second
fundamental form $A$, and $dv$ is the volume form for the induced
metric on $M$. Note that when $N=\mathbb{H}^{n+1}$,
$\overline{\Ric}(\nu,\nu)$ is equal to $-n$.

 In Section 3, we give some conditions for complete minimal hypersurfaces in
$\mathbb{H}^{n+1}$ to be stable as follows. If the $L^\infty$-norm
of the second fundamental form is sufficiently small at every
point in a complete minimal hypersurface $M$, then $M$ is stable
(Theorem \ref{thm:sufficient condition 1}). Moreover if the
$L^n$-norm of the second fundamental form is sufficiently small,
then $M$ is stable (Theorem \ref{thm:sufficient condition 2}).

In 1981, Mori \cite{Mori} explicitly described a one-parameter
family of complete stable minimal rotation surfaces in
$\mathbb{H}^3$. This example shows that a theorem due to do Carmo
and Peng \cite{dCP} and Fischer-Colbrie and Schoen \cite{FCS}
which says that a complete stable minimal surface in
$\mathbb{R}^3$ must be a plane, does not hold in $\mathbb{H}^3$.
Two years later, do Carmo and Dajczer \cite{dCD} found a larger
family of complete minimal rotation surfaces which are also
stable. In Section 4, we study stability of catenoids in the
hyperbolic space. In \cite{dCD}, it was shown that there exist a
one-parameter family of unstable catenoids $M_a$ in $\mathbb{H}^3$
for $1/2 < a< 0.69$. We improve the upper bound of $a$ by
estimating the $L^2$-norm of $|\nabla |A||$ in terms of the
$L^2$-norm of the second fundamental form $A$ (Theorem \ref{thm :
grad of |A|}). We also prove that the above unstable catenoid in
$\mathbb{H}^3$ should have index one (Theorem \ref{thm : index
1}). Recall that for a compact subset $\Omega$ in a complete
minimal hypersurface $M$ in $\mathbb{H}^{n+1}$, the {\it index} of
$\Omega$ is defined to be the number of negative eigenvalues of
the stability operator $L:=\Delta - |A|^2 + n$ on $\Omega$,
counting the multiplicity. The index of $M$ is defined as the
infimum of $\ind (\Omega)$ for all compact subset $\Omega$.
Moreover we provide a family of complete minimal hypersurfaces in
$\mathbb{H}^{n+1}$, which is an extension of Mori's result to
higher dimensional cases (Theorem \ref{thm : a family of
hyperbolic catenoids}). Finally we investigate stability of
helicoids in $\mathbb{H}^3$ in Section 5.
\section{First eigenvalue estimates}
In this section we first extend Candel's result to a simply
connceted complete minimal surface in a Riemannian manifold. The
proof is actually based on Candel's proof.
\begin{thm} \label{thm:manifold}
Let $\Sigma$ be a simply connected stable minimal surface in a
$3$-dimensional simply connected Riemannian manifold $N^3$ with
sectional curvature $K_N$ satisfying $-b^2 \leq K_N \leq -a^2 < 0$
for $0<a \leq b$. Then the first eigenvalue of $\Sigma$ satisfies
\begin{eqnarray*}
\frac{1}{4}a^2 \leq \lambda_1 (\Sigma) \leq \frac{4}{3}b^2.
\end{eqnarray*}
\end{thm}
\begin{pf}
First we find an upper bound for $\lambda_1(\Sigma)$. Let $\{e_1,
e_2, e_3 \}$ be orthonormal frames in $N$ such that the vectors
$\{e_1, e_2\}$ are tangent to $M$ and $e_3$ is normal to $M$. The
Gauss curvature equation implies that the sectional curvature
$K_\Sigma$ of $\Sigma$ satisfies
\begin{eqnarray} \label{ineq:curvature of Sigma}
K_\Sigma &=& R^1_{212} + h_{11}h_{22} - h_{12}^2 \nonumber \\
&=& R^1_{212} - \frac{|A|^2}{2} \leq -a^2 - \frac{|A|^2}{2} < 0,
\end{eqnarray}
where $R^1_{212}$ is the sectional curvature of $N$ for the
section determined by $e_1, e_2$ and $h_{ij}=\langle
\bar{\nabla}_{e_i}e_3,e_j\rangle$, $\bar{\nabla}$ denoting
Riemannian connection of $N$. Since $\Sigma$ is simply connected
and has negative Gaussian curvature, there are global polar
coordinates about any point in $\Sigma$. Using this polar
coordinates, the metric tensor $g$ of $\Sigma$ can be written as
\begin{eqnarray*}
g= dr^2 + \phi(r,\theta)^2 d\theta^2,
\end{eqnarray*}
where $\phi(0,\theta)=0$ and $\displaystyle{\frac{\partial
\phi}{\partial r} \Big|_{(0,\theta)} := \phi_r (0,\theta)=1}$.

Using the equality (\ref{ineq:curvature of Sigma}) and
$\overline{\Ric}(e_3) = R^3_{131} + R^3_{232}$, the stability
inequality (\ref{def : stability}) becomes
\begin{eqnarray} \label{ineq:stability under curvature assumption}
0 &\leq& \int_\Sigma |\nabla f|^2 -(|A|^2 + \overline{\Ric}(e_3))f^2 dv \nonumber\\
&\leq& \int_\Sigma |\nabla
f|^2 -(R^3_{131} + R^3_{232} +
2R^1_{212} - 2K_\Sigma)f^2 dv  \nonumber\\
&\leq& \int_\Sigma |\nabla f|^2 + 2K_\Sigma f^2 +4 b^2f^2 dv .
\end{eqnarray}
Since the inequality (\ref{ineq:stability under curvature
assumption}) holds for all compactly supported Lipschitz function
$f$ on $\Sigma$, we shall choose some specific functions which
depend only on the distance $r$ to the origin of the polar
coordinates in $\Sigma$. More precisely, given $R>0$, we consider
a family $\mathcal{F}$ of radial function $f$ such that $f(0)=0$,
$f(r)=0$ for $r\geq R>0$ and $f(r)$ is piecewise linear in $r$,
that is, $f''(r)=0$ except for finitely many values of $r$. Note
that the Gaussian curvature $\displaystyle{K_\Sigma =
-\frac{\phi_{rr}}{\phi}}$. Thus the inequality
(\ref{ineq:stability under curvature assumption}) can be written
as
\begin{eqnarray*}
2\int_0^{2\pi}\int_0^R f^2\phi_{rr} dr d\theta \leq
\int_0^{2\pi}\int_0^R f_r^2 \phi dr d\theta +4b^2
\int_0^{2\pi}\int_0^R f^2 \phi dr d\theta .
\end{eqnarray*}
Integrating the left side of the above inequality twice by parts
and using the properties of the function $f$, we obtain
\begin{eqnarray} \label{eqn:polar}
-\int_\Sigma K_\Sigma f^2 dv &=& \int_0^{2\pi}\int_0^R f^2
\phi_{rr} dr d\theta = \int_0^{2\pi}\Big[f^2 \phi_r \Big]_0^R d\theta - 2 \int_0^{2\pi}\int_0^R f
f_r \phi_r dr d\theta \nonumber \\
&=& \int_0^{2\pi}\Big[-2f f_r\phi\Big]_0^R d\theta + 2\int_0^{2\pi}\int_0^R (f f_r)_r \phi dr d\theta \nonumber \\
&=& 2\int_0^{2\pi}\int_0^R (f_r^2 \phi +f f_{rr}) \phi dr d\theta \nonumber \\
&=& 2\int_0^{2\pi}\int_0^R f_r^2
\phi dr d\theta = 2\int_\Sigma |\nabla f|^2 dv .
\end{eqnarray}
Combining the inequality (\ref{ineq:stability under curvature
assumption}) with the equation (\ref{eqn:polar}), we get
\begin{eqnarray*}
3\int_\Sigma |\nabla f|^2 dv \leq 4b^2 \int_\Sigma f^2 dv .
\end{eqnarray*}
Hence it follows that
\begin{eqnarray}  \label{ineq:upper bound}
\lambda_1 (\Sigma) \leq \inf_{f\in \mathcal{F}} \frac{\int_\Sigma
|\nabla f|^2 dv }{\int_\Sigma f^2 dv} \leq \frac{4}{3}b^2 .
\end{eqnarray}
Now we estimate a lower bound of $\lambda_1(\Sigma)$. The
Laplacian of the distance function $r$ on $\Sigma \subset N$
satisfies \cite{Choe}
\begin{eqnarray*}
\Delta r \geq a(2-|\nabla r|^2) \coth ar \geq a .
\end{eqnarray*}
Integrating both sides over $\Omega \subset \Sigma$, we get
\begin{eqnarray} \label{ineq:Cheeger}
a \area (\Omega) \leq \int_\Omega \Delta r dv = \int_{\partial
\Omega} \frac{\partial r}{\partial \nu} ds \leq \length (\partial
\Omega).
\end{eqnarray}
Recall that the {\it Cheeger} {\it constant} of a Riemannian manifold $M$, $h(M)$ is
defined by
\begin{eqnarray*}
h(M) := \inf_\Omega \frac{\length (\partial \Omega)}{\area
(\Omega)},
\end{eqnarray*}
where $\Omega$ ranges over all open submanifold of $M$, with
compact closure in $M$, and smooth boundary. Then applying
Cheeger's inequality \cite{Chavel} and inequality
(\ref{ineq:Cheeger}), we obtain
\begin{eqnarray} \label{ineq:lower bound}
\lambda_1 (\Sigma) \geq \frac{1}{4} h(\Sigma)^2 = \frac{1}{4}a^2 .
\end{eqnarray}
Therefore the theorem follows from (\ref{ineq:upper bound}) and
(\ref{ineq:lower bound}).  \qed
\end{pf}

The first eigenvalue of a complete minimal hypersurface in the
hyperbolic space is bounded below by a constant
$\frac{(n-1)^2}{4}$ as mentioned in the introduction. We give an
upper bound for a stable minimal hypersurface with finite
$L^2$-norm of the second fundamental form of $M$.

\begin{thm} \label{thm:first eigenvalue in H}
Let $M$ be a complete stable minimal hypersurface in
$\mathbb{H}^{n+1}$ with $\int_M |A|^2 dv  < \infty$. Then we have
\begin{eqnarray*}
\frac{(n-1)^2}{4} \leq  \lambda_1(M) \leq n^2.
\end{eqnarray*}
\end{thm}

\begin{Rmk}
There is no nontrivial example of such complete minimal
hypersurfaces in $\mathbb{R}^{n+1}$, since do Carmo and Peng \cite{dCP}
proved that a complete stable minimal hypersurface $M$ in
$\mathbb{R}^{n+1}$ with $\int_M |A|^2 dv  < \infty$ must be a
hyperplane. However, there exist several examples of complete
minimal hypersurfaces with finite $L^2$-norm of the second
fundamental form in the hyperbolic space as we will see in Section
4 and 5. Note that we do not assume that $M$ is simply connected,
which is different from Candel's result.
\end{Rmk}

\begin{pf}
It is sufficient to show that $\lambda_1 (M) \leq n^2$ by the
inequality (\ref{thm : CL}).

Take a function $f$ as follows. For a fixed point $p\in M$ and a
fixed $R>0$, define a function $f : M \rightarrow \mathbb{R}$ by
\begin{equation*}
f(x)=
 \begin{cases}
  1& \text{, $r(x)\leq R$,}\\
  2-\displaystyle{\frac{r(x)}{R}}& \text{, $R \leq r(x) \leq 3R$,}\\
  -1& \text{, $3R \leq r(x) \leq 4R$,}\\
  -5+\displaystyle{\frac{r(x)}{R}}& \text{, $4R \leq r(x) \leq 5R$,}\\
  0& \text{, $r(x) \geq 5R$,}
 \end{cases}
\end{equation*}
where $r(x)$ is the distance from $p$ to $x$ in $M$. Then it
follows that $\int_M f <0$. For $0 \leq t \leq R$, we now consider
a family of functions $\{f_t\}$ defined by
\begin{equation*}
f_t (x)=
 \begin{cases}
  1& \text{, $r(x)\leq R$,}\\
  2-\displaystyle{\frac{r(x)}{R}}& \text{, $R \leq r(x) \leq 2R +t$,}\\
  \displaystyle{-\frac{t}{R}}& \text{, $2R +t \leq r(x) \leq 4R+t$,}\\
  -5+\displaystyle{\frac{r(x)}{R}}& \text{, $4R+t \leq r(x) \leq 5R$,}\\
  0& \text{, $r(x) \geq 5R$.}
 \end{cases}
\end{equation*}
Then it is easy to see that there exists $t_0$, $0 < t_0 < R$,
such that $\int_M f_{t_0} = 0$. From the definition of $\lambda_1
(M)$ and $\lambda_1 (B_R)$ for a ball $B_R$ of radius $R$ centered
at $p$, it follows
\begin{eqnarray} \label{ineq : lambda_1}
\lambda_1 (M) \leq \lambda_1 (B_R) \leq \frac{\int_{B_R} |\nabla
\phi|^2}{\int_{B_R} \phi^2}
\end{eqnarray}
for any compactly supported Lipschitz function $\phi$ satisfying
$\int_{B_R} \phi = 0$.

 Now put $|A|f_{t_0}$ for $\phi$ in the inequality (\ref{ineq : lambda_1}). Then
\begin{eqnarray*}
\lambda_1(M) \int_{B_R} |A|^2 f_{t_0}^2 dv &\leq& \int_{B_R}
|\nabla
(|A|f_{t_0})|^2 dv \\
&=& \int_{B_R} |\nabla |A||^2 f_{t_0}^2 dv + \int_{B_R}
|A|^2|\nabla f_{t_0}|^2 dv +2\int_{B_R} |A| f_{t_0}
\langle\nabla|A|, \nabla f_{t_0}\rangle dv .
\end{eqnarray*}
Moreover, using Schwarz inequality, for any positive number
$\alpha>0$, we have
\begin{eqnarray*}
2\int_{B_R} |A| f_{t_0}  \langle\nabla|A|, \nabla f_{t_0}\rangle
dv \leq \alpha \int_{B_R} |A|^2|\nabla f_{t_0}|^2 dv +
\frac{1}{\alpha}\int_{B_R} |\nabla|A||^2 f_{t_0}^2 dv .
\end{eqnarray*}
Therefore we obtain
\begin{equation}  \label{ineq : test function 1}
\lambda_1(M) \int_{B_R} |A|^2 f_{t_0}^2 dv \leq
(1+\frac{1}{\alpha})\int_{B_R} |\nabla |A||^2 f_{t_0}^2 dv
+(1+\alpha) \int_{B_R} |A|^2|\nabla f_{t_0}|^2 dv .
\end{equation}
On the other hand,  Chern, do Carmo, and Kobayashi \cite{CdCK}
showed that
\begin{eqnarray} \label{eqn : CdCK}
\sum h_{i j} \Delta h_{i j} = -\sum h_{i j}^2 h_{k l}^2 -n\sum h_{i
j}^2.
\end{eqnarray}
Furthermore, we have
\begin{eqnarray} \label{eqn : laplacian of A}
|A|\Delta |A| + |\nabla |A||^2 = \frac{1}{2}\Delta |A|^2 = \sum
h_{i j k}^2 +  \sum h_{i j} \Delta h_{i j}
\end{eqnarray}
Combining (\ref{eqn : CdCK}) with (\ref{eqn : laplacian of A}), we
get
\begin{eqnarray*}
|A|\Delta |A| + |A|^4 + n|A|^2 =  |\nabla A|^2 -  |\nabla |A||^2.
\end{eqnarray*}
However the curvature estimate by Xin \cite{Xin} says that
\begin{eqnarray*}
|\nabla A|^2 - |\nabla |A||^2 \geq \frac{2}{n} |\nabla |A||^2,
\end{eqnarray*}
and hence we have
\begin{eqnarray*}
|A|\Delta |A| + |A|^4 + n|A|^2 \geq \frac{2}{n} |\nabla |A||^2.
\end{eqnarray*}
Multiplying both sides by a Lipschitz function $f^2$ with compact
support in $B_R \subset M$ and integrating over ${B_R}$, we have
\begin{eqnarray*}
\int_{B_R} f^2 |A|\Delta |A| dv + \int_{B_R} f^2 |A|^4 dv + n\int_{B_R} f^2
|A|^2 dv \geq \frac{2}{n} \int_{B_R} f^2 |\nabla |A||^2 dv
\end{eqnarray*}
The divergence theorem yields that
\begin{eqnarray*}
0 = \int_{B_R} \Div (|A|f^2\nabla |A|)dv = \int_{B_R} f^2 |A|\Delta |A|dv
+ \int_{B_R} |\nabla |A||^2 f^2 dv + 2\int_{B_R} |A| f \langle\nabla|A|,
\nabla f\rangle dv.
\end{eqnarray*}
Therefore
\begin{eqnarray} \label{ineq : f and A 01}
\int_{B_R} f^2 |A|^4 dv + n\int_{B_R} f^2 |A|^2 dv - \int_{B_R} |\nabla |A||^2
f^2 dv -2\int_{B_R} |A| f \langle\nabla|A|, \nabla f\rangle dv \nonumber \\
\geq \frac{2}{n} \int_{B_R} f^2 |\nabla |A||^2 dv.
\end{eqnarray}
Since $M$ is stable, we have
\begin{eqnarray*}
\int_M |\nabla \phi|^2 - (|A|^2 -n)\phi^2 dv \geq 0
\end{eqnarray*}
for any compactly supported function $\phi$ on $M$.
Substituting $|A|f$ for $\phi$ gives
\begin{eqnarray*}
\int_{B_R} |\nabla (|A|f)|^2 - (|A|^2 -n)|A|^2f^2 dv\geq 0.
\end{eqnarray*}
Thus
\begin{eqnarray} \label{ineq : f and A 02}
\int_{B_R} |A|^2|\nabla f|^2 dv+ \int_{B_R} |\nabla |A||^2f^2 dv + 2\int_{B_R}
|A| f \langle\nabla|A|, \nabla f\rangle dv \nonumber \\ \geq
\int_{B_R} |A|^4f^2 dv - n\int_{B_R} |A|^2 f^2 dv.
\end{eqnarray}
By (\ref{ineq : f and A 01}) and (\ref{ineq : f and A 02}), we
obtain, for any compactly supported Lipschitz function $f$
\begin{eqnarray} \label{ineq:stability}
\int_{B_R} |A|^2|\nabla f|^2 dv + 2n \int_{B_R} |A|^2f^2 dv \geq
\frac{2}{n}\int_{B_R} |\nabla |A||^2f^2 dv.
\end{eqnarray}
Combining (\ref{ineq : test function 1}) with the inequality
obtained by substituting $f_{t_0}$ for $f$ in
(\ref{ineq:stability}), we get
\begin{equation} \label{ineq : final}
\Big\{ 1 + \frac{2n(1+\alpha)}{\lambda_1(M)}\Big\} \int_{B_R}
|A|^2|\nabla f_{t_0}|^2 dv \geq \Big\{\frac{2}{n} -
\frac{2n(1+\frac{1}{\alpha})}{\lambda_1(M)}\Big\} \int_{B_R}
|\nabla |A||^2 f_{t_0}^2 dv .
\end{equation}
Now suppose that $\lambda_1 (M) > n^2$. Choosing $\alpha>0$
sufficiently large and letting $R\rightarrow \infty$, we obtain
$\nabla |A| \equiv 0$, i.e., $|A|$ is constant. However, since
$\int_M |A|^2 < \infty$ and the volume of $M$ is infinite, it
follows from the above inequality (\ref{ineq : final}) that
$|A|\equiv 0$ which means that $M$ is a totally geodesic
hyperplane. Since the first eigenvalue of totally geodesic
hyperplane in $\mathbb{H}^{n+1}$ is equal to
$\displaystyle{\frac{(n-1)^2}{4}}$, this is a contradiction.
Therefore we get $\lambda_1 (M) \leq n^2$. \qed
\end{pf}

\section{Sufficient conditions for stability of minimal hypersurfaces in $\mathbb{H}^{n+1}$}
In this section we prove that if $|A|$ is bounded by a
sufficiently small constant at every point in a complete minimal
hypersurface $M$ in the hyperbolic space, then $M$ must be stable.
More precisely,
\begin{thm} \label{thm:sufficient condition 1}
Let $M$ be a complete minimal hypersurface in $\mathbb{H}^{n+1}$.
If $|A| \leq \frac{(n+1)^2}{4}$ at every point in $M$, then $M$ is
stable.
\end{thm}
\begin{pf}
Since the lower bound of the first eigenvalue of $M$ is
$\frac{(n-1)^2}{4}$ by the inequality (\ref{thm : CL}), we have
\begin{eqnarray*}
\frac{(n-1)^2}{4} \leq \lambda_1 (M) \leq \frac{\int_M |\nabla
f|^2}{\int_M f^2}
\end{eqnarray*}
for every compactly supported Lipschitz function $f$ on $M$. Hence
the assumption that $|A|^2 \leq \frac{(n+1)^2}{4}$ implies
\begin{eqnarray*}
\int_M |\nabla f|^2 - (|A|^2 - n)f^2 dv \geq \int_M (\lambda_1 (M) +
n - |A|^2 )f^2 dv \geq 0,
\end{eqnarray*}
which completes the proof. \qed
\end{pf}

It is well-known that the following Sobolev inequality \cite{HS}
on a minimal hypersurface $M$ in $\mathbb{H}^{n+1}$ holds
\begin{eqnarray} \label{ineq:Sobolev}
\Big(\int_M |f|^{\frac{2n}{n-2}} dv\Big)^{\frac{n-2}{n}} \leq C_s
\int_M |\nabla f|^2 dv,
\end{eqnarray}
where $C_s$ is the Sobolev constant which dependent only on $n \geq 3$.
Using this inequality one obtains another sufficient condition for
minimal hypersurfaces to be stable.
\begin{thm} \label{thm:sufficient condition 2}
Let $M$ be a complete minimal hypersurface in $\mathbb{H}^{n+1}$, $n \geq 3$.
If $\int_M |A|^n dv \leq (\frac{1}{C_s})^{\frac{n}{2}}$, then $M$
is stable.
\end{thm}
\begin{pf}
It suffices to show that
\begin{eqnarray*}
\int_M |\nabla f|^2 - (|A|^2 - n)f^2 dv \geq 0
\end{eqnarray*}
for all compactly supported Lipschitz function $f$. By Sobolev
inequality (\ref{ineq:Sobolev}), we have
\begin{equation} \label{ineq:1}
\int_M |\nabla f|^2 - (|A|^2 - n)f^2 dv \geq \frac{1}{C_s}
\Big(\int_M |f|^{\frac{2n}{n-2}} dv\Big)^{\frac{n-2}{n}} - \int_M
|A|^2 f^2 dv .
\end{equation}
On the other hand, applying H\"{o}lder inequality, we get
\begin{eqnarray} \label{ineq:2}
\int_M |A|^2 f^2 dv \leq \Big(\int_M |A|^n dv \Big)^{\frac{2}{n}}
\Big(\int_M |f|^{\frac{2n}{n-2}} dv \Big)^{\frac{n-2}{n}}.
\end{eqnarray}
Combining (\ref{ineq:1}) with (\ref{ineq:2}) we have
\begin{eqnarray*}
\int_M |\nabla f|^2  - (|A|^2 - n)f^2 dv &\geq&
\Big\{\frac{1}{C_s} - \Big(\int_M |A|^n
dv\Big)^{\frac{2}{n}}\Big\} \Big(\int_M |f|^{\frac{2n}{n-2}} dv
\Big)^{\frac{n-2}{n}} \\
 &\geq& 0,   \hspace{1cm} \ \mbox{(by assumption)}
\end{eqnarray*}
which completes the proof. \qed
\end{pf}

\section{Catenoids in $\mathbb{H}^{n+1}$}
In \cite{dCD}, do Carmo and Dajczer proved that there exist three
types of rotationally symmetric minimal hypersurfaces in
$\mathbb{H}^{n+1}$. Following \cite{dCD}, we say that a
rotationally symmetric minimal hypersurface $M$ is a {\it
spherical} catenoid, if $M$ is foliated by spheres, a {\it
hyperbolic} catenoid, if it is foliated by hyperbolic spaces, and
a {\it parabolic} catenoid, if it is foliated by horospheres. Do
Carmo and Dajczer showed that the complete hyperbolic and
parabolic catenoids in $\mathbb{H}^3$ are all globally stable.
Furthermore they also proved that there exist some unstable
spherical catenoids in $\mathbb{H}^3$. In what follows, we shall
denote by $\mathbb{L}^{n+1}$ the space of $(n+1)$-tuples $x=(x_1,
\cdots, x_{n+1})$ with Lorentzian metric $\langle x, y\rangle =
-x_1y_1 +x_2y_2 + \cdots +x_{n+1}y_{n+1}$, where $y=(y_1, \cdots,
y_{n+1})$. The hyperbolic space $\mathbb{H}^n$ is the simply
connected hypersurface of $\mathbb{L}^{n+1}$ defined by
$\mathbb{H}^n = \{x\in \mathbb{L}^{n+1} : \langle x, x\rangle =
-1, x_1 \geq 1 \}$.

To state their result for unstable spherical catenoids in
$\mathbb{H}^3$, we parametrize a spherical catenoid in
$\mathbb{H}^3$ as follows. (See \cite{dCD} and \cite{Mori}.) For
each constant $a> 1/2$, define the mapping $f_a : \mathbb{R}
\times \mathbb{S}^1 \rightarrow \mathbb{H}^3$ by
\begin{eqnarray*}
f_a (s, \theta) = \Big( \sqrt{a \cosh(2s) + \frac{1}{2}} \cosh
\phi (s), \sqrt{a \cosh(2s) + \frac{1}{2}} \sinh \phi (s),\\
 \sqrt{a \cosh(2s) - \frac{1}{2}} \cos \theta, \sqrt{a \cosh(2s) - \frac{1}{2}} \sin \theta \Big),
\end{eqnarray*}
where $\displaystyle{\phi (s) = (a^2- \frac{1}{4})^{1/2}\int_0 ^s
\frac{1}{(a \cosh(2t) + \frac{1}{2})(a \cosh(2t) -
\frac{1}{2})^{1/2}} dt}$.\\

Do Carmo and Dajczer observed that if $1/2 < a < c_0, c_0 \simeq
0.69$, then the spherical catenoids $M_a$'s are unstable. We shall
improve the upper bound $c_0$ by using the inequality
(\ref{ineq:stability}), which is different from their method.
Letting $R \rightarrow \infty$ in (\ref{ineq:stability}), one can
immediately obtain the following.
\begin{thm} \label{thm : grad of |A|}
Let $M$ be a complete stable minimal hypersurface in
$\mathbb{H}^{n+1}$ with $\int_M |A|^2 dv  < \infty$. Then we have
\begin{eqnarray*}
\int_M |\nabla |A||^2 dv \leq n^2\int_M |A|^2 dv ,
\end{eqnarray*}
and hence the $L^2$-norm of $|\nabla |A||$ is finite.
\end{thm}

As a consequence of this Theorem \ref{thm : grad of |A|}, the
upper bound $c_0$ due to do Carmo and Dajczer can be improved as
follows.
\begin{cor} \label{cor : unstable}
Spherical catenoid $M_a$ in $\mathbb{H}^3$ is unstable for $1/2 <
a < c_0, \ c_0 \simeq 0.73$ .
\end{cor}

\begin{pf}
We first observe that the spherical catenoid $M_a$ satisfies
$\int_{M_a} |A|^2 dv < \infty$. To see this, we note that for
$a>1/2$
\begin{eqnarray*}
I &=& ds^2 + (a\cosh 2s - \frac{1}{2})dt^2, \\
|A|^2 &=& \frac{2(a^2-\frac{1}{4})}{(a \cosh 2s
-\frac{1}{2})^2}, \\
dv &=& (a \cosh 2s -\frac{1}{2})^{\frac{1}{2}}dsdt \ \ \text{for $a>\frac{1}{2}$ and $0 \leq t \leq 2\pi$.}\\
\end{eqnarray*}
Thus
\begin{eqnarray*}
\int_M |A|^2 dv &=& 8\pi (a^2-\frac{1}{4})\int_0^\infty
\frac{1}{(a
\cosh 2s -\frac{1}{2})^{\frac{3}{2}}} ds  \\
&<& 8\pi (a^2-\frac{1}{4})\int_0^\infty \frac{1}{(a +a
s^2-\frac{1}{2})^{\frac{3}{2}}} ds  < \infty.
\end{eqnarray*}
Now define a function $F(a)$ by
\begin{eqnarray*}
F(a) : = 4 \int_M |A|^2 dv - \int_M |\nabla |A||^2 dv.
\end{eqnarray*}
Using $\displaystyle{|\nabla |A|| = \Big|-\sqrt{2(a^2
-\frac{1}{4})}\frac{2a\sinh 2s}{(a\cosh 2s
-\frac{1}{2})^2}\Big|}$, we have
\begin{eqnarray*}
F(a) = 32\pi (a^2 - \frac{1}{4}) \int_0^\infty \Big\{ \frac{1}{(a
\cosh 2s -\frac{1}{2})^{\frac{3}{2}}} - a^2 \frac{\sinh^2
2s}{(a\cosh2s -\frac{1}{2})^{\frac{7}{2}}}\Big\}ds.
\end{eqnarray*}
By Theorem \ref{thm : grad of |A|}, we see that if $M_a$ is stable
for some $a$, then $F(a) \geq 0$. However a straightforward
computation by using a computer shows that $F(a)
 < 0$ for $1/2 < a < c_0$, $c_0 \simeq 0.73$. Therefore we get the
 conclusion. \qed
\end{pf}

As we have seen before, there exist some unstable catenoids in
$\mathbb{H}^3$. Hence it is interesting to find the index of such
catenoids which measures the degree of instability. It is
well-known that catenoids have index $1$ in $\mathbb{R}^3$. Very
recently, Tam and Zhou \cite{TZ} proved that higher dimensional
catenoids in $\mathbb{R}^{n+1}$ with $n \geq3$ have index one.
Motivated by this, we shall prove the following result using the
similar arguments as in \cite{TZ}.

\begin{thm} \label{thm : index 1}
Let $M$ be a spherical catenoid in $\mathbb{H}^{n+1}$. Then the
index of $M$ is at most $1$.
\end{thm}
\begin{pf}
We may assume that $M$ is unstable. It suffices to show that the
second eigenvalue $\lambda_2 (D) \geq 0$ of the stability operator
$L = \Delta + |A|^2 - n$ on some bounded domain $D\subset M$. We
prove this theorem by contradiction. For this purpose, suppose
 that the index of $M$ is at least $2$. Then there exists a domain $D(R)=(-R,R) \times
\mathbb{S}^{n-1}$ such that $\lambda_2 (D(R)) < 0$ for $R>0$.

Let $f$ be the second eigenfunction satisfying
\begin{eqnarray*}
Lf&=&-\lambda_2 (D(R)) f \hspace{2cm} {\rm in \ }\ \ D(R)\\
f&=&0 \hspace{3.5cm} {\rm on \ }  \ \ \partial D(R) .
\end{eqnarray*}
We claim that $f$ is rotationally symmetric, that is, $f(t_1,
\cdots, t_{n-1}, s) = f(s)$. \\
To see this, consider a generating curve $\alpha (s) := (x(s),
y(s), z(s), 0, \cdots, 0) \subset \mathbb{H}^{n+1}$ and its
rotation axis $\{(\cosh u, \sinh u, 0, \cdots, 0)\}\subset
\mathbb{H}^{n+1}$. Let $P_0$ be the totally geodesic hyperplane
such that $P_0 \perp \alpha'(0)$ and $\alpha (0) = (1, 0, \cdots,
0) \in P_0$. For any vector $v\in S_{\alpha(0)} P_0 := \{v \in
T_{\alpha(0)} P_0 : |v|=1\}$, denote by $P_v$ the (unique) totally
geodesic hyperplane such that $\alpha(0)\in P_v$ and $P_v \perp v$
at $\alpha(0)$.

Let $\sigma_v$ be the reflection across the hyperplane $P_v$. For
any point $p\in D(R)$, define the difference function $\varphi_v
(t_1, \cdots, t_{n-1}, s)$ by
\begin{eqnarray*}
\varphi_v (t_1, \cdots, t_{n-1}, s) := f(t_1, \cdots, t_{n-1}, s)
- f_v(t_1, \cdots, t_{n-1}, s),
\end{eqnarray*}
where $f_v (p) := f(\sigma_v (p))$. Then it follows that $\Delta f
= \Delta f_v$. Thus
\begin{equation} \label{eqn : lambda_2}
 \begin{cases}
  L \varphi_v &= -\lambda_2 (D(R)) \hspace{2cm} {\rm in \ }\ \ D(R)\\
  \varphi_v &= 0 \hspace{3.5cm} {\rm on \ }  \ \ \partial D(R)\cap
P_v .\\
 \end{cases}
\end{equation}

Since $P_v$ divides $D(R)$ into two parts, we choose one of them
and denote by $D_v^+(R)$. Note that $D_v^+(R)$ is a minimal graph
over a domain $P_v$. Hence $D_v^+(R)$ is stable. However from
(\ref{eqn : lambda_2}) and the assumption that $\lambda_2<0$, it
follows that $\varphi_v\equiv 0$. As in the Euclidean space, any
rotaion around the axis $\{(\cosh u, \sinh u, 0, \cdots, 0)\}
\subset \mathbb{H}^{n+1}$ can be expressed as a composition of
finite number of reflections. Since $v$ was arbitrarily chosen,
the claim is obtained.

On the other hand, since the second eigenfunction of the operator
$L$ changes sign, there exists a number $r_0 \in (-R,R)$
satisfying $f(r_0)=0$. We may assume that $r_0 \geq 0$ and the
second eigenfunction $f(s)> 0$ on the domain $D(r_0, R) = \{(t_1,
\cdots, t_{n-1}, s) \in D(R) : s\in (r_o, R)\}$. The function $f$
is still an eigenfunction of $L$ on $D(r_0, R)$. Moreover it is
easy to see that $D(r_0, R)$ is a minimal graph over the
hyperplane $P_0$, which means that $D(r_0, R)$ is stable. This is
a contradiction to the assumption that $\lambda_2 <0$. Therefore
we get the conclusion. \qed
\end{pf}

\begin{Rmk}
When $n=2$, we observed that a spherical catenoid $M_a \subset
\mathbb{H}^3$ is unstable if $1/2 < a < 0.73$ in Corollary
\ref{cor : unstable}. It follows from the above theorem that these
spherical catenoids must have index 1.
\end{Rmk}

We now describe stability of hyperbolic catenoids in the
hyperbolic space $\mathbb{H}^{n+1}$. For that purpose, we give a
parametrization of a hyperbolic catenoid generated by a curve
$(x(s),y(s),z(s))$ in the hyperbolic plane $\mathbb{H}^2$ which is
parametrized by arclength. It follows that
\begin{eqnarray}
-x(s)^2 + y(s)^2 + z(s)^2 &=& -1 , \hspace{1cm} x(s) \geq 1  \label{eqn:included in H}\\
-x'(s)^2 + y'(s)^2 + z'(s)^2 &=& 1  \label{eqn:arclength},
\end{eqnarray}
\begin{eqnarray}
f(t_1, \cdots, t_{n-1}, s) = (x(s)\varphi_1, \cdots,
x(s)\varphi_n, y(s), z(s)), \label{hyperbolic catenoid}\\
\varphi_i = \varphi_i (t_1, \cdots, t_{n-1}), \ \ -\varphi_1^2 +
\varphi_2^2+ \cdots + \varphi_n^2 = -1, \nonumber
\end{eqnarray}
where $(\varphi_1, \cdots, \varphi_n)$ is an orthogonal
parametrization of the hyperbolic space $\mathbb{H}^{n-1}$. From
(\ref{eqn:included in H}) and (\ref{eqn:arclength}), $y(s)$ and
$z(s)$ are determined by
\begin{eqnarray*}
y(s) &=& \sqrt{x(s)^2-1}\sin \phi(s),\\
z(s) &=& \sqrt{x(s)^2-1}\cos \phi(s),
\end{eqnarray*}
where $\displaystyle{\phi(s)= \int_0^s
\frac{\sqrt{x^2-x'^2-1}}{x^2-1}dt}$.

Using minimality and rotationally symmetric property of a
catenoid, one can see that the direction of the parameters are
principal directions and the principal curvatures are given by
\begin{eqnarray*}
\lambda_1 = \cdots = \lambda_{n-1} =
-\frac{\sqrt{x^2-x'^2-1}}{x},\\
\lambda_n = \frac{x''-x}{\sqrt{x^2-x'^2-1}} =
(n-1)\frac{\sqrt{x^2-x'^2-1}}{x}.
\end{eqnarray*}
(See \cite[Proposition 3.2]{dCD}.) Furthermore we can write down
an ordinary differential equation as follows \cite[Lemma
3.15]{dCD}:
\begin{eqnarray}\label{eqn:ode}
x' = \sqrt{x^2-1-a^2 x^{2(1-n)}}, \hspace{2cm} a=\mbox{const}.
\end{eqnarray}
To find a unique solution of (\ref{eqn:ode}), we fix initial data
as follows:
\begin{eqnarray*}
x(0) &=& t \ \ \geq 1 \\
x'(0) &=& 0.
\end{eqnarray*}
Then from the initial data it follows that
\begin{eqnarray} \label{eqn : constant a or t}
a = t^{n-1}\sqrt{t^2-1} \geq 0 .
\end{eqnarray}
Moreover in order to have a nontrivial parametrization of a
hyperbolic catenoid we see that $a>0$. Therefore for each constant
$t>1$ the parametrization $f(t_1, \cdots, t_{n-1}, s)$ defines a
hyperbolic catenoid $M_t$ in $\mathbb{H}^{n+1}$. We now state a
our result about stability of hyperbolic catenoids in the
hyperbolic space.
\begin{thm} \label{thm : a family of hyperbolic catenoids}
Let $M_t$ be a family of hyperbolic catenoids in
$\mathbb{H}^{n+1}$ defined as in (\ref{hyperbolic catenoid}). Then
$M_t$ is a complete stable hypersurface in $\mathbb{H}^{n+1}$ for
$1<t<1+\frac{(n+1)^2}{4n(n-1)}$.
\end{thm}
\begin{pf}
Observe that
\begin{eqnarray*}
|A|^2 &=& \sum \lambda_i ^2 = (n-1)\lambda_1 ^2 + \lambda_n ^2 = (n-1) \lambda_1 ^2 + (n-1)^2 \lambda_1 ^2\\
&=& n(n-1)\lambda_1 ^2 \\
&=& n(n-1)\frac{x^2-x'^2-1}{x^2}\\
&=& n(n-1)\frac{a^2}{x^{2n}} \hspace{27mm}\mbox{(by (\ref{eqn:ode}))} \\
&=& n(n-1) \frac{t^{2(n-1)}(t^2-1)}{x^{2n}} . \hspace{1cm}
\mbox{(by (\ref{eqn : constant a or t}))}
\end{eqnarray*}
Since $x(s)$ is monotonically increasing by (\ref{eqn:ode}), we
get $x(s) \geq x(0) = t > 1$. Therefore $|A|^2 \leq
n(n-1)(t^2-1)$. The assumption on $t$ implies that $|A|^2 \leq
\frac{(n+1)^2}{4}$. The conclusion follows from Theorem
\ref{thm:sufficient condition 1}. \qed
\end{pf}

\begin{Rmk}
It is not hard to see that the family $\{M_t\}$ of hyperbolic catenoids in the above theorem satisfy $\int_{M_t} |A|^2 dv < \infty$ by using 
$\sqrt{x^2-1-a^2} < x' < \sqrt{x^2-1}$, which is obtained from equality (\ref{eqn:ode}) and the fact that $x>1$. 
\end{Rmk}
\section{Helicoids in $\mathbb{H}^3$}
Let $l$ be a geodesic in $\mathbb{H}^3$. Let $\{\psi_t\}$ be the
translation of distance $t$ along $l$ and let $\{\varphi_t\}$ be
the rotation of angle $t$ around $l$. Given any $\alpha \in
\mathbb{R}$, one can see that $\lambda = \{\lambda_t\} = \{\psi_t
 \circ \varphi_{\alpha t}\}$ is a one-parameter subgroup of
 isometries of $\mathbb{H}^3$ which is called a {\it helicoidal} group
 of isometries with angular pitch $\alpha$. A {\it helicoid} in
 $\mathbb{H}^3$ is a $\lambda$-invariant surface. (See
 \cite{Ripoll}.) In 1989, Ripoll \cite{Ripoll} proved that a
 helicoid $M_\alpha$ with angular pitch $|\alpha|<1$ is stable by
 showing that such $M_\alpha$ foliates $\mathbb{H}^3$. In this
 section, we improve the upper bound of angular pitch $|\alpha|$ by
 simple arguments.

A helicoid $M_\alpha \subset \mathbb{H}^3 \subset \mathbb{L}^4$
can be written explicitly as follows \cite{BDJ} :
\begin{eqnarray*}
X(s,t) = (\cosh s \cosh t, \sinh s \cosh t, \cos\alpha s \sinh t,
\sin \alpha s \sinh t).
\end{eqnarray*}
A little computation shows that the first and second fundamental
forms of $M_\alpha$ are given by
\begin{eqnarray*}
I &=& (\cosh^2 t + \alpha^2 \sinh^2 t)ds^2 + dt^2,\\
II &=& -2\frac{\alpha}{\sqrt{\cosh^2 t + \alpha^2 \sinh^2 t}} \
dsdt .
\end{eqnarray*}
Since $\cosh^2 t + \alpha^2 \sinh^2 t \geq 1$, it follows
\begin{eqnarray*}
|A|^2 = \frac{\alpha^2}{\cosh^2 t + \alpha^2 \sinh^2 t} +
\frac{\alpha^2}{(\cosh^2 t + \alpha^2 \sinh^2 t)^3} \leq
2\alpha^2.
\end{eqnarray*}
The following theorem is an immediate consequence of Theorem
\ref{thm:sufficient condition 1}.
\begin{thm}
A helicoid $M_\alpha$ with angular pitch $|\alpha|^2 \leq
\frac{9}{8}$ is stable.
\end{thm}

\end{document}